\documentclass[12pt,a4paper]{amsart}
\usepackage[utf8]{inputenc}	
\usepackage[T1]{fontenc}
\usepackage[in,headings]{fullpage}
\usepackage{amsbsy, amscd, amsfonts, amsmath, amsrefs, amssymb, amsthm}
\usepackage{graphicx}
\usepackage{float}
\usepackage{color}   
\usepackage[colorlinks=true,linkcolor=blue,citecolor=blue,urlcolor=blue]{hyperref}

\newtheorem{theorem}{Theorem}

\newtheorem{corollary}[theorem]{Corollary}
\newtheorem{proposition}[theorem]{Proposition}
\newtheorem{lemma}[theorem]{Lemma}

\theoremstyle{definition}

\newtheorem{remark}[theorem]{Remark}

\theoremstyle{remark}

\newcommand{\C}{\mathbf{C}}
\newcommand{\Z}{\mathbf{Z}}

\newcommand{\R}{\mathbf{R}}
\newcommand{\N}{\mathbf{N}}

\renewcommand{\Re}{\mathop{\mathrm{Re}}\nolimits}

\newcommand{\Rzeta}{\mathop{\mathcal R }\nolimits}

\newfont{\cmbsy}{cmbsy10}
\newfont{\cmmib}{cmmib10}
\newcommand{\Orden}{\mathop{\hbox{\cmbsy O}}\nolimits}

\begin{document}

\title{Density theorems for Riemann's auxiliary function. }
\author[Arias de Reyna]{J. Arias de Reyna}
\address{%
Universidad de Sevilla \\ 
Facultad de Matem\'aticas \\ 
c/Tarfia, sn \\ 
41012-Sevilla \\ 
Spain.} 

\subjclass[2020]{Primary 11M99; Secondary 30D99}

\keywords{zeta function, Riemann's auxiliar function}


\email{arias@us.es, ariasdereyna1947@gmail.com}


\begin{abstract}
We prove a density theorem for the auxiliar function $\Rzeta(s)$ found by Siegel in Riemann papers. Let $\alpha$ be a real number with $\frac12< \alpha\le 1$, and let $N(\alpha,T)$ be the number of zeros $\rho=\beta+i\gamma$ of $\Rzeta(s)$ with $1\ge \beta\ge\alpha$ and $0<\gamma\le T$. Then we prove \[N(\alpha,T)\ll T^{\frac32-\alpha}(\log T)^3.\]

Therefore, most of the zeros of $\Rzeta(s)$ are near the critical line or to the left of that line. The imaginary line for $\pi^{-s/2}\Gamma(s/2)\Rzeta(s)$ passing through a zero of $\Rzeta(s)$ near the critical line frequently will cut the critical line, producing two zeros of $\zeta(s)$ in the critical line. 
\end{abstract}

\maketitle
\section{Introduction}

Siegel \cite{Siegel} introduced the function $\Rzeta(s)$ considered by Riemann in his Nachlass.
It is connected with $\zeta(s)$ by the relation 
\[Z(t)=e^{i\vartheta(t)}\zeta(\tfrac12+it)=2\Re\{e^{i\vartheta(t)}\Rzeta (\tfrac12+it)\}.\]
Therefore, $\zeta(s)$ have a zero on the critical line at each point $s=\frac12+it$  where $e^{i\vartheta(t)}\Rzeta (\tfrac12+it)$  is purely imaginary.  This is the same as saying that
$\zeta(s)$ have a zero at a point $s$ in the critical line just in case $\pi^{-s/2}\Gamma(s/2)\Rzeta(s)$ is purely imaginary.  

Siegel connected the zeros of $\Rzeta(s)$ with zeros on the critical line of $\zeta(s)$, showing that $\zeta(s)$ have at least $2N_2(T)$ zeros on the critical line with ordinates between $0$ and $T$, where $N_2(T)$ is the number of zeros of $\Rzeta(s)$ to the left of the critical line and with ordinates between $0$ and $T$ . To understand the relationship between the zeros, consider the x-ray of 
$e^{i\vartheta(t)}\Rzeta(\frac12+it)$ in the figure. The horizontal line is the real axis for $t$.
\begin{figure}[H]
\begin{center}
\includegraphics[width=0.9\hsize]{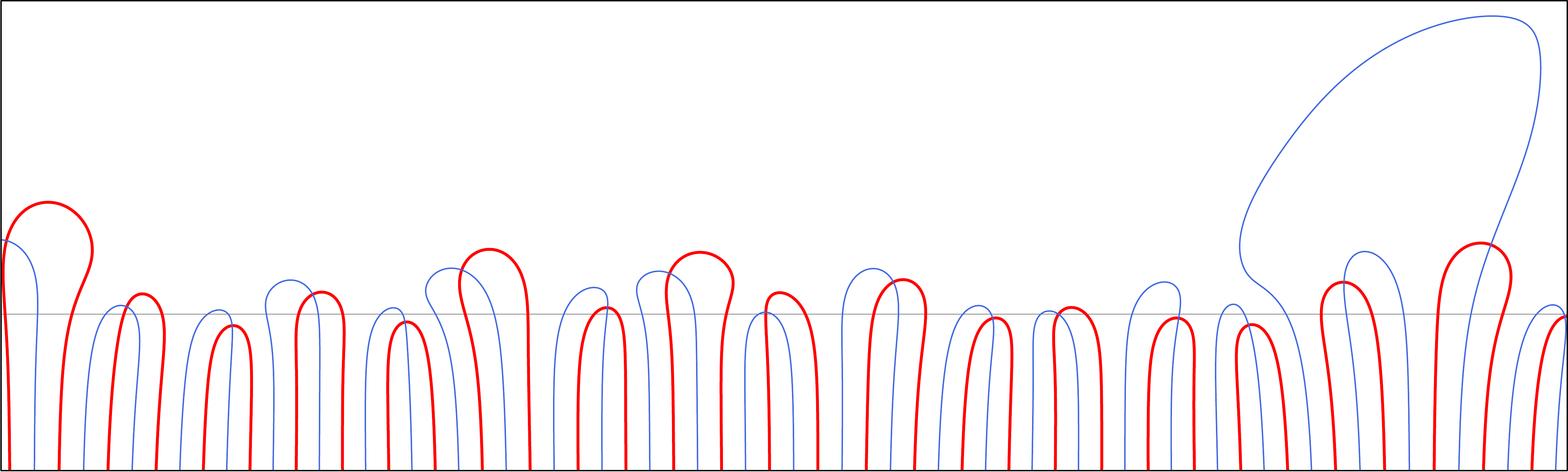}
\caption{x-ray of $e^{i\vartheta(t)}\Rzeta(\frac12+it)$ in $(200\,040,200\,060)\times(-2,4)$}
\label{default}
\end{center}
\end{figure}
Each zero of $\Rzeta(s)$ to the left of the critical line gives a zero above the real axis in the figure. The function $Z(t)$ vanishes at the points where the imaginary lines (blue lines) in the figure cut the real axis. 
All imaginary lines (in blue) start at the bottom of the figure so that the imaginary line passing through a zero of $e^{i\vartheta(t)}\Rzeta(\frac12+it)$
above the real axis cuts the real axis at two points, producing two zeros of $Z(t)$. 

In principle, it would be possible that the imaginary line through a zero below the real axis in the figure do not cut the real axis. But we see in the examples in the figure that they are sufficiently close to give also two zeros of $Z(t)$.

Since $\Rzeta(s)$ can be approximated by  Dirichlet polynomials, the methods used to prove the density theorems for $\zeta(s)$ can be applied to $\Rzeta(s)$, and this is the objective of this paper. Our main result is Theorem \ref{T:firstthm} where we prove that for $1/2<\alpha\le1$, we have $N(\alpha,T)\ll T^{\frac32-\alpha}\log^3T$, where $N(\alpha,T)$ denotes the number of zeros of $\Rzeta(s)$ in $[\alpha,1]\times[0,T]$ counted with multiplicity.  The proof is standard, sure this can be improved, I have not used a mollifier here. 

These zeros below the real line in our figure correspond to the zeros of $\Rzeta(s)$ to the right of the critical line. In this paper, we show that most of the zeros of $\Rzeta(s)$ to the right of the critical line are  near this line. 
These density theorems  make it almost inevitable that each of these zeros produce its two zeros for $\zeta(s)$. 
In this way, the $\frac{T}{4\pi}\log\frac{T}{2\pi}-\frac{T}{4\pi}$ zeros of $\Rzeta(s)$ can generate the $\frac{T}{2\pi}\log\frac{T}{2\pi}-\frac{T}{2\pi}$ zeros  predicted by the Riemann hypothesis.

It is also interesting that the zeros of  $\Rzeta(s)$ give limits to the possible density theorems.

In this paper, we follow the Theorems in Iwaniec \cite{Iw},
adapting them to  $\Rzeta(s)$. I give proofs of these Theorems giving explicit constants. But this is not needed for our Theorems. 

We use the notations $N(\sigma,T)$, that usually refers to the zeta function, for the corresponding functions associated with $\Rzeta(s)$. Since we do not refer to the corresponding function for $\zeta(s)$ I expect this to not cause confusion. We denote by $\N=\{1,2,\dots\}$ the set of natural numbers.

\section{Dirichlet Polynomials and Detection of Zeros}

We follow the Theorems in Iwaniec \cite{Iw} but give his theorems with definite constants. This is not necessary. An expert can proceed directly to the Section \ref{S:density}. The only notable difference with the case of the zeta function is that we cannot directly use the interval $[0,T]$ because the Dirichlet polynomials only approximates $\Rzeta(s)$ in intervals of the type $[2\pi N^2,2\pi(N+1)^2]$ or close to this interval.

First, we give two elementary lemmas on extension of differentiable functions.

\begin{lemma}\label{P:190805-1}
Given three real numbers $a$, $b$ and $c$, with  $|a|\le1$, $|b|\le1$ and $|c-b|\le1$, there is a function  $\varphi\colon[0,1]\to\R$ of class $\mathcal C^2$ such that
\[\varphi(0)=\varphi'(0)=\varphi''(0)=0,\quad \varphi(1)=a,\quad \varphi'(1)=b,\quad \varphi''(1)=c,\]
satisfying the following inequalities
\begin{equation}\label{E:190805-2}
|\varphi(x)|\le \tfrac{10}{9},\quad |\varphi'(x)|\le \tfrac73,\quad |\varphi''(x)-\varphi'(x)|\le \tfrac{41}{5},\qquad 0\le x\le 1.\end{equation}
\end{lemma}

\begin{proof}
Let $d=c-b$, then we want a function $\varphi\in\mathcal C^2[0,1]$ such that
\begin{equation}\label{E:200510-1}
\varphi(0)=\varphi'(0)=\varphi''(0)=0,\quad \varphi(1)=a,\quad \varphi'(1)=b,\quad
\varphi''(1)-\varphi'(1)=d,
\end{equation}
and with the smaller possible values of the three norms $\Vert\varphi\Vert_\infty$, 
$\Vert\varphi'\Vert_\infty$ and $\Vert\varphi''-\varphi'\Vert_\infty$
Here $-1\le a, b, d\le 1$ are arbitrary. Any $(a,b,d)$ is a convex combination of 
points $(\varepsilon_1,\varepsilon_2,\varepsilon_3)$ with $\varepsilon_j=\pm\,1$. Therefore, we search functions $\varphi_j\in\mathcal C^2[0,1]$ for $1\le j\le 4$ verifying 
\begin{align*} 
&\varphi_1(0)=\varphi_1'(0)=\varphi_1''(0)=0, \quad \varphi_1(1)=1;\quad \varphi'_1(1)=1;\quad
\varphi''_1(1)-\varphi'_1(1)=1.\\
&\varphi_2(0)=\varphi_2'(0)=\varphi_2''(0)=0, \quad \varphi_2(1)=1;\quad \varphi'_2(1)=1;\quad
\varphi''_2(1)-\varphi'_2(1)=-1.\\
&\varphi_3(0)=\varphi_3'(0)=\varphi_3''(0)=0, \quad \varphi_3(1)=1;\quad \varphi'_3(1)=-1;\quad
\varphi''_3(1)-\varphi'_3(1)=1.\\
&\varphi_4(0)=\varphi_4'(0)=\varphi_4''(0)=0, \quad \varphi_4(1)=1;\quad \varphi'_4(1)=-1;\quad
\varphi''_4(1)-\varphi'_4(1)=-1.
\end{align*}
According to what we have said, given $-1\le a, b, d\le 1$ there are eight real numbers $a_j\ge0$ with
$\sum_{j=1}^8 a_j=1$ and such that $\varphi=\sum_{j=1}^4 a_j\varphi_j-\sum_{j=1}^4 a_{j+4}\varphi_j$ satisfies the conditions \eqref{E:200510-1}.
And we obtain our Lemma taking the supremum of the norms of each $\varphi_j$ as bounds.
We  take 
\begin{align*}
\varphi_1(x)&=x^3(4x^2-10x+7), & \varphi_3(x)&=\frac{x^3}{5}(21x^6-28x^5-28x^4+178x^2-250x+112),\\
\varphi_2(x)&=x^3(3x^2-8x+6), &  \varphi_4(x)&=\frac{x^3}{10}(49x^6-70x^5-56x^4+381x^2-522x+228).
\end{align*}
and we obtain 
\begin{align*}
&\Vert\varphi_1\Vert_\infty=1, &&\Vert\varphi_1'\Vert_\infty=1.52797, &&\Vert\varphi_1''-\varphi_1'\Vert_\infty=3.68937,\\
&\Vert\varphi_2\Vert_\infty=1, &&\Vert\varphi_2'\Vert_\infty=\frac{189}{125}, &&\Vert\varphi_2''-\varphi_2'\Vert_\infty=3.35572,\\
&\Vert\varphi_3\Vert_\infty=1.10624, &&\Vert\varphi_3'\Vert_\infty=2.29980, &&\Vert\varphi_3''-\varphi_3'\Vert_\infty=8.18258,\\
&\Vert\varphi_4\Vert_\infty=1.09127, &&\Vert\varphi_4'\Vert_\infty=2.23453, &&\Vert\varphi_4''-\varphi_4'\Vert_\infty=8.11817.\qedhere\\
\end{align*}
\end{proof}

\begin{lemma}\label{L:2}
Let $f\colon[1,N]\to\R$ of class $\mathcal C^2$, such that
\[|f(x)|\le1,\quad x|f'(x)|\le 1,\quad x^2 |f''(x)|\le 1,\qquad 1\le x\le N.\]
There exists a function $g\colon(0,+\infty)\to\R$ of class $\mathcal C^2$, with support contained in  $[1/e,eN]$, which extends  $f$, and such that
\[|g(x)|\le \tfrac{10}{9},\quad x|g'(x)|\le \tfrac73,\quad x^2 |g''(x)|\le \tfrac{41}{5},\qquad 0< x< +\infty.\]
\end{lemma}

\begin{proof}
Define $F$ on $[0,\log N]$ by $F(\log x)=f(x)$.  We have
\[f(x)=F(\log x), \quad f'(x)=\frac{1}{x} F'(\log x),\quad f''(x)=\frac{1}{x^2}(F''(\log x)-F'(\log x))\]
Therefore, for  $y$ in the domain of $F$ we have
\[|F(y)|\le 1,\quad |F'(y)|\le 1,\quad |F''(y)-F'(y)|\le 1.\]
These inequalities for $F$ and its derivatives are equivalent to the hypothesis for $f$.

To prove the lemma, we construct an extension  of $F$,   $G\colon\R\to\R$ of class $\mathcal C^2$, with support contained in  $[-1, \log N+1]$, and such that
\[|G(y)|\le \tfrac{10}{9},\quad |G'(y)|\le \tfrac73,\quad |G''(y)-G'(y)|\le \tfrac{41}{5}.\]
It suffices to construct a function on $[-1,0]$ such that its derivatives to the order $2$ coincide with the ones of $F$ at $y=0$, and are $0$ at the point $-1$. And an analogous construction of an adequate function on $[\log N, \log N+1]$.

To define $G$ we extends $F$ to $[-1, \log N+1]$ with these values and then $G(y)=0$ for
$y\not\in[-1, \log N+1]$. 

To construct $G$ on $[-1,0]$ we apply
Proposition \ref{P:190805-1} with the  constants $a=F(0)$, $b=F'(0)$, 
$c=F''(0)$. These values satisfy by hypothesis the inequalities $|a|\le1$, $|b|\le1$ y $|c-b|\le1$.

Proposition \ref{P:190805-1}  gives us a function $\varphi\colon[0,1]\to\R$. Define then $G$ en $[-1,0]$ by means of the expression
\[G(x)=\varphi(1+x),\quad  G'(x)=\varphi'(1+x),\quad G''(x)=\varphi''(1+x).\]
We will have
\[G(0)=\varphi(1)=F(0), \quad G'(0)=\varphi'(1)=F'(0),\quad
G''(0)=\varphi''(1)=F''(0),\]
and
\[G(-1)=G'(-1)=G''(-1)=0.\]
According to Proposition \ref{P:190805-1} for  $-1\le x\le 0$ we have
\[|G(x)|\le \tfrac{10}{9}, \quad |G'(x)|\le\tfrac73, \quad |G''(x)-G'(x)|\le \tfrac{41}{5}.\]
so that $g(x)=G(\log x)$ satisfies the conditions in our Lemma. 
\end{proof}

Using this lemma, Mellin's transform, and integration by parts, we arrive at the following bound 
of a Dirichlet polynomial.
\begin{proposition}\label{P:Iwaniec1}
Let  $N\ge2$ be a natural number and  $f\colon[1,N]\to\R$ a $\mathcal C^2$-function  
satisfying the inequalities
\[|f(x)|\le 1,\quad x|f'(x)|\le 1,\qquad x^2|f''(x)|\le 1.\]
Then
\[\Bigl|\sum_{n=1}^N a_n f(n)n^{-it_0}\Bigr|\le \frac{b}{2\pi}(\log N+2) \int_{-\infty}^{+\infty}\Bigl|\sum_{n=1}^N
a_n n^{-i(t_0+t)}\Bigr|\frac{dt}{(1+t^2)},\qquad b<8.775.\]
\end{proposition}

\begin{proof}
By Lemma \ref{L:2} we may extend $f$ to a $\mathcal C^2$ class  function $g\colon(0,+\infty)\to\R$ with the support contained in $[1/e,eN]$ and such that 
\[|g(x)|\le \tfrac{10}{9},\quad x|g'(x)|\le \tfrac73,\qquad x^2|g''(x)|\le \tfrac{41}{5}.\]
Consider its Mellin transform
\[h(t)=\int_0^\infty g(x)x^{it}\frac{dx}{x}.\]
Let $M=41/5$ and $m=10/9$ be the constants in the lemma \ref{L:2}.
There is an algebraic number $a=2.626198\dots$ such that $M/m=a\sqrt{1+a^2}$, so that 
\[b:=\frac{M}{a}\sqrt{1+a^2}=m(1+a^2), \qquad b=8.77435\dots\]

Integrating by parts 
\begin{multline*}
h(t)=\Bigl. g(x)\frac{x^{it}}{it}\Bigr|_0^\infty-\frac{1}{it}\int_0^\infty g'(x) x^{it}\,dx=
-\frac{1}{it}\int_0^\infty g'(x) x^{it}\,dx\\=\Bigl. -\frac{1}{it}g'(x)\frac{x^{it+1}}{1+it}\Bigr|_0^\infty+\frac{1}{it(1+it)}\int_0^\infty g''(x)x^{it+1}\,dx.
\end{multline*}
It follows that 
\[|h(t)|\le \frac{1}{|t(1+it)|}\int_0^\infty x^2|g''(x)|\frac{dx}{x}\le 
\frac{M}{|t(1+it)|}\int_{1/e}^{eN}\frac{dx}{x}=\frac{M(\log N+2)}{|t(1+it)|}.\]
Then for $|t|>a$
\[|h(t)|\le \frac{M(\log N+2)}{(1+t^2)}\frac{1+t^2}{|t(1+it)|}\le \frac{M(1+t^2)}{|t|\sqrt{1+t^2}}\frac{(\log N+2)}{1+t^2}\le \frac{b(\log N+2)}{1+t^2}.\]

For  $|t|\le a$   we have, from its definition, the following
\[|h(x)|\le \int_{1/e}^{eN}|g(x)|\frac{dx}{x}\le m\log(N+2)\le m(1+a^2)\frac{(\log N+2)}{1+t^2}= b\frac{(\log N+2)}{1+t^2}.\]

Once estimated  $h(x)$, we finish the proof using  Mellin's inversion
\[g(x)=\frac{1}{2\pi}\int_{-\infty}^{+\infty} h(t)x^{-it}\,dt.\]
Since  $f(n)=g(n)$ for $1\le n\le N$, we obtain
\[\Bigl|\sum_{n=1}^N a_n f(n)n^{-it_0}\Bigr|=\Bigl|\sum_{n=1}^N  \frac{a_n n^{-it_0}}{2\pi}\int_{-\infty}^{+\infty} h(t)n^{-it}\,dt\Bigr|=\Bigl|  \frac{1}{2\pi}\int_{-\infty}^{+\infty} h(t)\sum_{n=1}^Na_n n^{-it_0-it}\,dt\Bigr|\le\]
\[\le \frac{1}{2\pi}\int_{-\infty}^{+\infty}\Bigl|\sum_{n=1}^Na_n n^{-i(t_0+t)}\Bigr|\cdot|h(t)|\,dt
\le \frac{b(\log N+2)}{2\pi}\int_{-\infty}^{+\infty}\Bigl|\sum_{n=1}^Na_n n^{-i(t_0+t)}\Bigr|\frac{dt}{1+t^2}.\qedhere
\]
\end{proof}

Usually it is easier to use the next corollary.
\begin{corollary}\label{C:Iwaniec2}
Under the same conditions as in Proposition \ref{P:Iwaniec1} we have
\[\Bigl|\sum_{n=1}^N a_n f(n)n^{-it_0}\Bigr|\le \frac{b}{2\sqrt{\pi}}(\log N+2) \Bigl(\int_{-\infty}^{+\infty}\Bigl|\sum_{n=1}^N
a_n n^{-i(t_0+t)}\Bigr|^2\frac{dt}{1+t^2}\Bigr)^{1/2}.\]
\end{corollary}
\begin{proof}
Applying Schwarz's Lemma
\[\int_{-\infty}^{+\infty} u(t)v(t)\frac{dt}{1+t^2}\le \Bigl(\int_{-\infty}^{+\infty} u(t)^2\frac{dt}{1+t^2}\Bigr)^{1/2}\Bigl(\int_{-\infty}^{+\infty} v(t)^2\frac{dt}{1+t^2}\Bigr)^{1/2}.\qedhere\]
\end{proof}

To apply this corollary we need a Theorem estimating the mean of a Dirichlet polynomial.
\begin{theorem}[Montgomery]
Suppose that $\lambda_1$, \dots, $\lambda_N$ are distinct real numbers, for any complex numbers $\{x_n\}_{n=1}^N$ and $\{y_n\}_{n=1}^N$ we have 
\begin{equation}
\Bigl|\sum_{\substack{1\le m\le N\\1\le n\le N\\ n\ne m}}\frac{x_m y_n}{\lambda_m-\lambda_n}\Bigr|\le \frac{3\pi}{2}\Bigl(\sum_{m=1}^N |x_m|^2/\delta_m\Bigr)^{1/2}\Bigl(\sum_{n=1}^N |y_n|^2/\delta_n\Bigr)^{1/2},
\end{equation}
where
$\displaystyle{\delta_n=\min_{\substack{1\le m\le N\\m\ne n}}|\lambda_m-\lambda_n|}$.
\end{theorem}
See Montgomery \cite{Mont}*{p.~140, eq.~(27)}. The constant $3\pi/2$ is not the best possible. Preissman \cite{P} proved that the inequality is true with the constant 
$\pi\sqrt{1+\frac23\sqrt{\frac65}}<\frac43\pi$.

\begin{theorem}\label{T:Montgomery}
For any real numbers $U$, $T>0$ and complex numbers $a_n$, we have
\begin{equation}
\int_U^{U+T} \Bigl|\sum_{n=1}^N \frac{a_n}{n^{it}}\Bigr|^2\,dt\le (T+\tfrac{4\pi}{3})\sum_{n=1}^N |a_n|^2+ \tfrac{8\pi}{3}\sum_{n=1}^N n|a_n|^2.
\end{equation}
\end{theorem}

\begin{proof}
We may assume that $U=0$ by changing $a_n$ into $a_n n^{-iU}$. 
We have
\[\int_0^T \Bigl|\sum_{n=1}^N \frac{a_n}{n^{it}}\Bigr|^2\,dt=\sum_{n,m=1}^N
a_n\overline{a_m}\int_0^T e^{it(\log m-\log n)}\,dt\]\[=
T\sum_{n=1}^N |a_n|^2+\sum_{n\ne m}\frac{a_n \overline{a_m}}{i(\log m-\log n)}(e^{iT(\log m-\log n)}-1).\]
Applying two times Montgomery's Theorem \ref{T:Montgomery}
(With $\lambda_n=\log n$, we have  $\delta_1=\log 2$ and   $\delta_n=\log \frac{n+1}{n}\ge \frac{1}{n+1/2}$ for $n\ge1$)
\[\le (T+\tfrac{4\pi}{3})\sum_{n=1}^N |a_n|^2+2\frac{4\pi}{3} \sum_{n=1}^N n|a_n|^2.\qedhere\]
\end{proof}

All the above results combine to get a bound on the set of points where a Dirichlet polynomial gets large values. This will be our main tool to obtain estimates on the number of zeros.

\begin{lemma}\label{L:simpleineq}
Let $T>2$ and $0\le t_1< t_2<\cdots<t_n\le T$, with $t_{k+1}-t_k\ge1$ for $1\le k\le n-1$, then we have 
\[\sum_{k=1}^n \frac{1}{1+(t-t_k)^2}\le\pi\coth(\pi), \qquad\text{ for all $t\in\R$}.\]
and for $mT\le t\le (m+1)T$ or $-mT\le t\le -(m-1)T$ with $m\ge 2$ we have
\[\sum_{k=1}^n \frac{1}{1+(t-t_k)^2}\le \frac{1}{(m-1)^2T^2}+\arctan((m-1)T)-\arctan((m-2)T).\]
\end{lemma}
\begin{proof}
In any case, there are points $t'_k\in[0,T]$ with $t'_{k+1}=t'_k+1$ and such that 
$|t-t'_k|\le |t-t_k|$, then,  with $x=t-t'_1+1$, we have
\[S:=\sum_{k=1}^n \frac{1}{1+(t-t_k)^2}\le \sum_{k=1}^n \frac{1}{1+(t-t'_k)^2}=
\sum_{k=1}^n \frac{1}{1+(t-t'_1-k+1)^2}\]\[\le \sum_{k\in\Z}\frac{1}{1+(t-t'_1-k+1)^2}
=\frac{\pi\sinh(2\pi)}{\cosh(2\pi)-\cos(2\pi x)}\\ \le \frac{\pi\sinh(2\pi)}{\cosh(2\pi)-1}=\pi\coth(\pi)\]

The case $-mT\le t\le -(m-1)T$ is entirely similar to the on in which $mT\le t\le (m+1)T$. Hence we only consider the cases  $mT\le t\le (m+1)T$ with $m\ge2$ an integer.
In the case $mT\le t\le (m+1)T$ we have with $x=t-t'_1+1$
\[S\le \sum_{k=1}^n \frac{1}{1+(t-t'_k)^2}=\sum_{k=1}^n \frac{1}{1+(x-k)^2}\]
We have 
\[(m-1)T+1=mT-T+1\le x=t-t'_1+1\le (m+1)T-0+1.\]
Hence for $m\ge 2$ we have $x>T+1$ and $\frac{1}{1+(x-u)^2}$ is an increasing function of $u$ so that 
\[S\le\sum_{k=1}^n \frac{1}{1+(x-k)^2}\le \frac{1}{1+(x-n)^2}+\int_1^n \frac{du}{1+(x-u)^2}\]\[=\frac{1}{1+(x-1)^2}+\arctan(x-1)-\arctan(x-n)\]
Now $(m-1)T\le x-1\le (m+1)T$, and $T\ge t_n-0>n-1$, so that  $x-n\ge (m-1)T+1-(T+1)=(m-2)T\ge0$. 
Therefore,
\[S\le \frac{1}{1+(m-1)^2T^2}+\arctan((m+1)T)-\arctan((m-2)T).\qedhere\]
\end{proof}

\begin{theorem}\label{T:Iwaniec4}
Let  $A(s)=\sum_{n=1}^N \frac{a_n}{n^s}$ be a  Dirichlet polynomial.  Let $s_1$, $s_2$, \dots, $s_J$ be a finite set of well separated  points, i.e.  points $s_j=\sigma_j+it_j$ such that $|t_j-t_{j'}|\ge1$ for $j\ne j'$. Assume also that the points are in the rectangle $R=[0,\Delta]\times[0,T]$ with $\Delta+\Delta^2\le 1$ and $T\ge4$.

Then 
\begin{equation}
\sum_{j=1}^J |A(s_j)|^2
\le25\pi (\log N+2)^2\Bigl\{(T+\tfrac{4\pi}{3})\sum_{n=1}^N |a_n|^2+ \tfrac{8\pi}{3}\sum_{n=1}^N n|a_n|^2\Bigr\}.
\end{equation}
\end{theorem}

\begin{proof}
Write $A(s_j)$ as 
\[A(s_j)=\sum_{n=1}^N \frac{1}{n^{\sigma_j}} a_n n^{-i t_j}=
\sum_{n=1}^N f_j(n) a_n n^{-i t_j},\]
where $f_j(x):=x^{-\sigma_j}$. 
For $x\ge1$ we have
\[|f_j(x)|\le 1,\quad f'_j(x)=-\frac{\sigma_j}{x}f_j(x),\quad f_j''(x)=\Bigl(\frac{\sigma_j}{x^2}+
\frac{\sigma_j^2}{x^2}\Bigr)f_j(x),\]
Since $\Delta^2+\Delta\le 1$, the functions $f_j$ satisfy the conditions of Proposition \ref{P:Iwaniec1} and Corollary \ref{C:Iwaniec2}. Let $A(s)=\sum_{n=1}^N a_n n^{-s}$, then 
\begin{multline*}\frac{4\pi}{b^2(\log N+2)^2}\sum_{j=1}^J|A(s_j)|^2\\\le \sum_{j=1}^J \int_{-\infty}^\infty |A(it_j+it)|^2\frac{dt}{1+t^2}= \int_{-\infty}^\infty |A(it)|^2\Bigl(\sum_{j=1}^J \frac{1}{1+(t-t_j)^2}\Bigr)\,dt
\end{multline*}
We decompose the integral in subintervals of length $T$
\[\int_{-\infty}^\infty |A(it)|^2\Bigl(\sum_{j=1}^J \frac{1}{1+(t-t_j)^2}\Bigr)\,dt=
\sum_{m=-\infty}^{\infty}\int_{mT}^{(m+1)T} |A(it)|^2\Bigl(\sum_{j=1}^J \frac{1}{1+(t-t_j)^2}\Bigr)\,dt.\]
We apply the lemma \ref{L:simpleineq}. For the intervals $[-T,0]$, $[0,T]$ and $[T,2T]$ we use the bound $\pi\coth(\pi)$. For the other, use the alternative inequality

Therefore,
\begin{align*}
P&:=\int_{-\infty}^\infty |A(it)|^2\Bigl(\sum_{j=1}^J \frac{1}{1+(t-t_j)^2}\Bigr)\,dt\\
&\le \pi\coth\pi\int_{-T}^0 |A(it)|^2\,dt+
\pi\coth\pi\int_0^T |A(it)|^2\,dt+\pi\coth\pi\int_T^{2T} |A(it)|^2\,dt
\end{align*}
\begin{align*}
&
+\sum_{m=2}^\infty\Bigl(\frac{1}{1+(m-1)^2T^2}+\arctan((m-1)T)-\arctan((m-2)T)\Bigr)\int_{kT}^{(k+1)T}|A(it)|^2\,dt\\
&+\sum_{m=2}^\infty\Bigl(\frac{1}{1+(m-1)^2T^2}+\arctan((m-1)T)-\arctan((m-2)T)\Bigr)\int_{-mT}^{-(m-1)T}\mskip-10mu|A(it)|^2\,dt.
\end{align*}
By Theorem \ref{T:Montgomery} all these integrals have the same bound. 
\[\int_{kT}^{(k+1)T}|A(it)|^2\le H:= (T+\tfrac{4\pi}{3})\sum_{n=1}^N |a_n|^2+ \tfrac{8\pi}{3}\sum_{n=1}^N n|a_n|^2, \qquad k\in\Z.\]

And we obtain
\[P\le 
H\Bigl(3\pi\coth\pi+2\sum_{m=2}^\infty\Bigl(\frac{1}{(m-1)^2T^2}+\arctan((m-1)T)-\arctan((m-2)T)\Bigr)\Bigr)\]

The other sum is bounded in the same way. Therefore,
\[P\le H\Bigl(3\pi\coth\pi+2\frac{\zeta(2)}{T^2}+\pi\Bigr).\]
This means that 
\[\sum_{j=1}^J|A(s_j)|^2\le \frac{b^2}{4\pi}\Bigl(3\pi\coth\pi+2\frac{\zeta(2)}{T^2}+\pi\Bigr)(\log N+2)^2 H.\]
For $T\ge4$
\[\frac{b^2}{4\pi}\Bigl(3\pi\coth\pi+2\frac{\zeta(2)}{T^2}+\pi\Bigr)=24.97621 \pi\le 25\pi.\qedhere\]
\end{proof}

\section{First density theorem}\label{S:density}

For $\sigma<1$ and $T>0$ we denote by $N(\sigma,T)$ the number of zeros of $\Rzeta(s)$ in the rectangle $[\sigma,1]\times[0,T]$ counted with its multiplicity. We know \cite{A173} that all zeros 
$\rho=\beta+i\gamma$ of $\Rzeta(s)$ with $\gamma>0$ and $\gamma> t_0$ satisfies $\beta<1$, so that $N(\sigma,T)$ counts essentially all zeros with $\gamma>0$ and $\beta\ge\sigma$. This is proved in \cite{A173} for an astronomical $t_0$, but I believe \cite{A172},  that this is also true taking $t_0=0$.  Since we do not speak about zeros of $\zeta(s)$ in this paper, the notation $N(\sigma,T)$ do not lead to any confussion.

Siegel proved that the zeros $\rho=\beta+i\gamma$ of $\Rzeta(s)$  with $\gamma>0$ satisfies $-\beta\le\gamma^{3/7}$ and asserts that he could prove that 
$-\beta\le\gamma^{\varepsilon}$ for any $\varepsilon>0$. I have not been able to prove the Siegel assertion. In \cite{A98}*{Theorem 7} I proved that there is no zero
in the closed set 
\[F:=\{s=\sigma+it\colon t>5408\pi, \   (1-\sigma)>\sqrt{225\pi t}\}.\]

We will need a bound for the number of zeros contained in a rectangle $[\alpha,1]\times[T,T+1]$ of height 1. 

\begin{theorem}\label{T:densityzeros}
Let $\alpha\le 1$ and $T\ge e$, there is an absolute constant $C$ such that the number of zeros $\rho=\beta+i\gamma$ of $\Rzeta(s)$  with $\alpha\le \beta\le1$ and $T\le \gamma\le T+1$, counted with multiplicity, is less than $C(3-\alpha)\log T$. 
\end{theorem}

\begin{proof} 
First assume that $T>T_0>20000\pi$ and let $N(T)$ be the number of zeros   with $\alpha\le \beta\le1$ and $T\le \gamma\le T+1$. For $T_1<T_2$, denote by $N(\alpha, T_1,T_2)=N(\alpha,T_2)-N(\alpha,T_1)$. If $(1-\alpha)>\sqrt{225\pi (T+1)}$, the result quoted  in \cite{A98}*{Theorem 7} implies that 
\[N(\alpha,T, T+1)=N(1-\sqrt{225\pi (T+1)},T,T+1).\]
Hence, we only need to prove our theorem for $(1-\alpha)\le \sqrt{225\pi (T+1)}$. And we assume this inequality and $T\ge T_0$  in what follows.

We apply Theorem D in the book by Ingham \cite{In}*{p.~49}. Consider a disc of center at $2+i(T+\frac12)$ passing through the point $\alpha+iT$. The radius is determined by $r^2=(2-\alpha)^2+1/4$ and depends only on $\alpha$. Let $D$ be the disc with the same center and radius $=2r$. It is clear that $D\subset R:=[2-2r,2+2r]\times[T+\frac12-2r, T+\frac12+2r]$. Part of this rectangle $R$ is contained in $\sigma\ge0$, where by \cite{A92}*{Prop.~12}  we have $|\Rzeta(s)|\le 2\sqrt{t/2\pi}$. In the rectangle with $\sigma\le 0$ we apply \cite{A92}*{Prop.~13}
\[|\Rzeta(\sigma+it)|\le
\frac{19\,t}{(2\pi)^{1-\sigma}}\{(1-\sigma)^2+t^2\}^{\frac{1}{4}-\frac{\sigma}{2}}\qquad
\sigma\le 0,\quad t\ge 16\pi.\]
To apply this, we need $T+\frac12-2r>16\pi$, but notice that for $T>6000\pi$
\[16\pi+2r<16\pi+4(2-\alpha)<16\pi+4+4\sqrt{225\pi(T+1)}<T+\tfrac12,\]
so that $T+\frac12>16\pi+2r$.
Therefore 
\[|\Rzeta(\sigma+it)|\le \frac{19(2\pi)^\sigma}{2\pi}\Bigl(1+\frac{(1-\sigma)^2}{t^2}\Bigr)^{\frac14-\frac{\sigma}{2}}t^{\frac32-\sigma}.\]
We have $1-\sigma<t$. To prove it, observe that for $s$ in the rectangle $R$ we have $t>T+\frac12-2r$ and $2-2r\le \sigma$. Therefore,
$1-\sigma\le 2r-1$ and we only have to prove $2r-1\le T+\frac12-2r$. Since 
$r<2(2-\alpha)$,  we have 
\[4r-\tfrac32\le 8(2-\alpha)- \tfrac32=\tfrac{13}{2}+8(1-\alpha)\le 
\tfrac{13}{2}+8\sqrt{225\pi (T+1)}.\]
Therefore $(1-\sigma)\le t$ follows from 
\[\tfrac{13}{2}+8\sqrt{225\pi (T+1)}<T.\]
That is true for $T>20000\pi$.

Hence, for $s\in D$ with $\sigma\le0$ we have
\begin{align*}
|\Rzeta(\sigma+it)|&\le\frac{19}{2\pi}2^{\frac14-\frac{\sigma}{2}} (T+\tfrac12+2r)^{\frac32-\sigma}
\le \frac{19}{2\pi}2^{\frac14-\frac{\sigma}{2}} \Bigl(1+\frac{\frac12+2r}{1-\sigma+2r-\frac12}\Bigr)^{\frac32-\sigma}
T^{\frac32-\sigma}\\
&\le \frac{19}{2\pi}2^{\frac74-\frac{3\sigma}{2}}T^{\frac32-\sigma}.\end{align*}
For $s\in D$, we have 
\[\frac32-\sigma\le 2r-\frac12\le 2(2-\alpha)+\frac12\le 2(3-\alpha).\]
Hence, for $s\in D$ with $\sigma<0$, 
\[|\Rzeta(\sigma+it)|\le\frac{19}{2\pi}2^{3(3-\alpha)}T^{2(3-\alpha)}.\]
The inequality is also true for $s\in D$ with $\sigma>0$, since $\alpha\le 1$ and $|\Rzeta(\sigma+it)|\le 2(t/2\pi)^{1/2}$ in this case. 

So, $\Rzeta(s)$ is bounded in $D$ by $c2^{3(3-\alpha)} T^{2(3-\alpha)}$. The number  $N$ of zeros of $\Rzeta(s)$ in the initial disc satisfies, by Theorem D in Ingham 
\[2^N\le \frac{c2^{3(3-\alpha)}}{|\Rzeta(2+i(T+\frac12))|}T^{2(3-\alpha)}.\]
By \cite{A173}*{Prop.~6} we know that for $T\ge T_0$ we have $|\Rzeta(2+i(T+\frac12))|\ge1/8$.  

Hence,
\[N\log2\le 2(3-\alpha)\log T+(3-\alpha)\log 8+\log(8c).\]
From which we derive 
\[N\log 2\le (3-\alpha)\log T\Bigl(2+\frac{\log 8}{\log T}+\frac{\log (8c)}{(3-\alpha)\log T}\Bigr)\le
(3-\alpha)\log T\Bigl(2+\frac{\log 8}{\log T}+\frac{\log (8c)}{2\log T}\Bigr) \]
For any $\varepsilon>0$ there is a $T_0(\varepsilon)$ such that 
\[N\le (2+\varepsilon)(3-\alpha)\frac{\log T}{\log 2},\qquad T\ge T'_0(\varepsilon).\]
 
Taking $\varepsilon=1$ for example, we know that there is a constant $C$ such that 
$N(T)\le C(3-\alpha)\log T$ for any $\sigma\le 1$ and $T\ge T_0$. The number of zeros of  $\Rzeta(s)$ with $e\le\gamma\le T_0$ is finite (by \cite{A98}*{Theorem 7}), so that $N(T)\le C(3-\alpha)\log T$ is still true in this case for $T>e$ taking $C$ large enough.
\end{proof}

In \cite{A173}*{Prop.~5} it is proved that for $0\le \sigma\le1$ and $t\ge3\pi$  we have
\[\Rzeta(s)=\sum_{n\le \sqrt{t/2\pi}}\frac{1}{n^s}+\Orden(t^{-\sigma/2}),\]
where the implicit constant is absolute.
Therefore, for $2\pi K^2\le t<2\pi(K+1)^2$,  where $K$ is a natural number, we have
\[\Rzeta(s)=\sum_{n=1}^K\frac{1}{n^s}+\Orden(K^{-\sigma}), \qquad K^2\le \frac{t}{2\pi}\le(K+1)^2.\]
We denote by $N_K(\sigma)$ the number of zeros $\rho=\beta+i\gamma$ of $\Rzeta(s)$ contained in the rectangle $[\sigma,1]\times(2\pi K^2, 2\pi(K+1)^2]$. Notice that by the results in \cite{A173}
for large $K$ it is the same as the number of zeros in the rectangle $[\sigma,+\infty)\times(2\pi K^2, 2\pi(K+1)^2]$.

\begin{theorem}\label{T:firstthm}
For $\frac12< \alpha\le 1$, let $N(\alpha,T)$ be the number of zeros $\rho=\beta+i\gamma$ of $\Rzeta(s)$ with $1\ge \beta\ge\alpha$ and $0<\gamma\le T$. Then 
\[N(\alpha,T)\ll T^{\frac32-\alpha}(\log T)^3.\]
\end{theorem}

\begin{proof}
First, consider a natural number $K$ and estimate the number of zeros $N_K(\alpha)$ contained in the rectangle $Q_K=[\alpha,1]\times(2\pi K^2, 2\pi(K+1)^2]$.  Notice that for $s\in Q_K$ we have
$t\sim K^2$ so that $\log t\sim \log K$.  Divide $Q_K$ into $\Orden(K)$ rectangles of height $1$, by Theorem \ref{T:densityzeros}, in each of these rectangles there are at most $c\log K$ zeros. Therefore, we may get  $J\gg N_K(\alpha)/\log K$ \emph{well separated} zeros of $\Rzeta(s)$, that is, such that the difference from two of its heights is always $\ge 1$. Let $\rho_j=\beta_j+i\gamma_j$, $1\le j\le J$ these zeros. 

Consider the Dirichlet polynomial 
\[D(s)=\sum_{n=2}^K \frac{1}{n^{\alpha+2\pi i K^2}}\frac{1}{n^s}\]
For each $j$ we have 
\[0=\Rzeta(\rho_j)=1+D(\rho_j-\alpha-2\pi i K^2)+\Orden(K^{-\alpha}).\]
Since the implicit  constant in $\Orden(K^{-\alpha})$ is absolute, there is a $K_0$ such that
for $K\ge K_0$ the above relation implies 
\[|D(\rho_j-\alpha-2\pi i K^2)|\ge\tfrac12.\]
Notice that the points $\rho_j-\alpha-2\pi i K^2=x_j+iy_j$ with $0\le x_j=\beta_j-\alpha\le 1-\alpha<\frac12$ and $0\le y_j=\gamma_j-2\pi K^2\le 2\pi(2K+1)$.  We are in conditions to apply
Theorem \ref{T:Iwaniec4}, obtaining
\[\frac{J}{4}\le 25\pi(\log K+2)^2\Bigl\{(2\pi(2K+1)+\tfrac{4\pi}{3})\sum_{n=2}^K\frac{1}{n^{2\alpha}}+\tfrac{8\pi}{3}\sum_{n=2}^K\frac{n}{n^{2\alpha}}
\Bigr\}\ll K(\log K)^2 K^{1-2\alpha}.\]
From which we derive 
\begin{equation}\label{E:210218-1}
N_K(\alpha)\ll  K^{2-2\alpha}(\log K)^3,\qquad K\ge K_0.
\end{equation}
Since  each $N_K(\alpha)$ is finite, this is also true  for any $K\ge 2$, perhaps with a bigger implicit constant.

Let now $T>8\pi$ be given. There is a natural number $K\ge2$ such that $2\pi K^2\le T<2\pi (K+1)^2$. 
Then 
\begin{align*}
N(\alpha, T)&\le N(\alpha,8\pi)+\sum_{k=2}^K (N(\alpha,2\pi(k+1)^2-N(\alpha,2\pi k^2))\\
&\le 
N(\alpha,8\pi)+\sum_{k=2}^K N_k(\alpha)\ll
N(\alpha,8\pi)+\sum_{k=2}^K k^{2-2\alpha}(\log k)^3\\
&\ll K^{3-2\alpha}(\log K)^3\ll T^{\frac32-\alpha}(\log T)^3.\qedhere
\end{align*}
\end{proof}

\begin{proposition}
Let $f\colon[0,+\infty)\to(0,+\infty)$ be a nonincreasing function with \[\lim_{t\to+\infty}f(t)=0\quad\text{and}\quad\frac{f(t)\log t}{\log\log t}>4, \quad t\ge t_0.\]
The number of zeros of $\Rzeta(s)$ contained in the subset
\[\Omega_T:=\bigl\{s=\sigma+it\in \C\colon e\le t\le T, \sigma\ge \tfrac12+f(t)\bigr\}\]
is $o(T)$. Therefore, most of the zeros of $\Rzeta(s)$ are on the left of the line
$\sigma=\frac12+f(t)$. 
\end{proposition}

\begin{proof}
Consider a value of $K$ such that $2\pi K^2<T\le2\pi(K+1)^2$. All zeros in 
$\Omega_T$ with $2\pi k^2<\gamma\le 2\pi(k+1)^2$ are contained in 
$Q_k=[\alpha_k,1]\times(2\pi k^2,2\pi(k+1)^2]$ with $\alpha_k=\frac12+f(2\pi(k+1)^2)$. Hence, the number of zeros $N(T)$ of $\Rzeta(s)$ in $\Omega_T$ is less than or equal to 
\[N(T)\le C+\sum_{k=2}^KN_k(\tfrac12+f(2\pi(k+1)^2))\]
and by \eqref{E:210218-1} we get  
\[N(T)\ll C+\sum_{k=2}^K k^{1-f(2\pi(k+1)^2)}(\log k)^3\le C+\sum_{k=2}^K k(\log k)^3 \exp\Bigl(-f(2\pi(k+1)^2)\log k\Bigr)\]
\begin{align*}
N(T)&\le C+\sum_{k=2}^K \frac{k(\log k)^3}{(\log k)^4}\exp\bigl(4\log\log k-f(2\pi(k+1)^2)\log k\bigr)\\ &\ll
\sum_{k=2}^K \frac{k}{(\log k)^{1-\varepsilon}}\ll\frac{K^2}{(\log K)^{1-\varepsilon}}\ll T/(\log T)^{1-\varepsilon}.\qedhere
\end{align*}
\end{proof}

\begin{remark}
The next figures are the xray of $\pi^{-s/2}\Gamma(s/2)\Rzeta(s)$ for $\sigma\in (-105,105)$ and $t\in (-105,105)$, $(105,315)$ and $(315,525)$. The horizontal lines are for $t=2\pi n^2$ for the first values of $n$. We have drawn also the critical line. Each cut of the critical line with an imaginary line is a zero of $\zeta(s)$. There are 140 zeros of $\Rzeta(s)$ represented in these xrays, and 287 zeros of zeta.
\end{remark}

\begin{figure}[H]
\begin{center}
\includegraphics[width=\hsize]{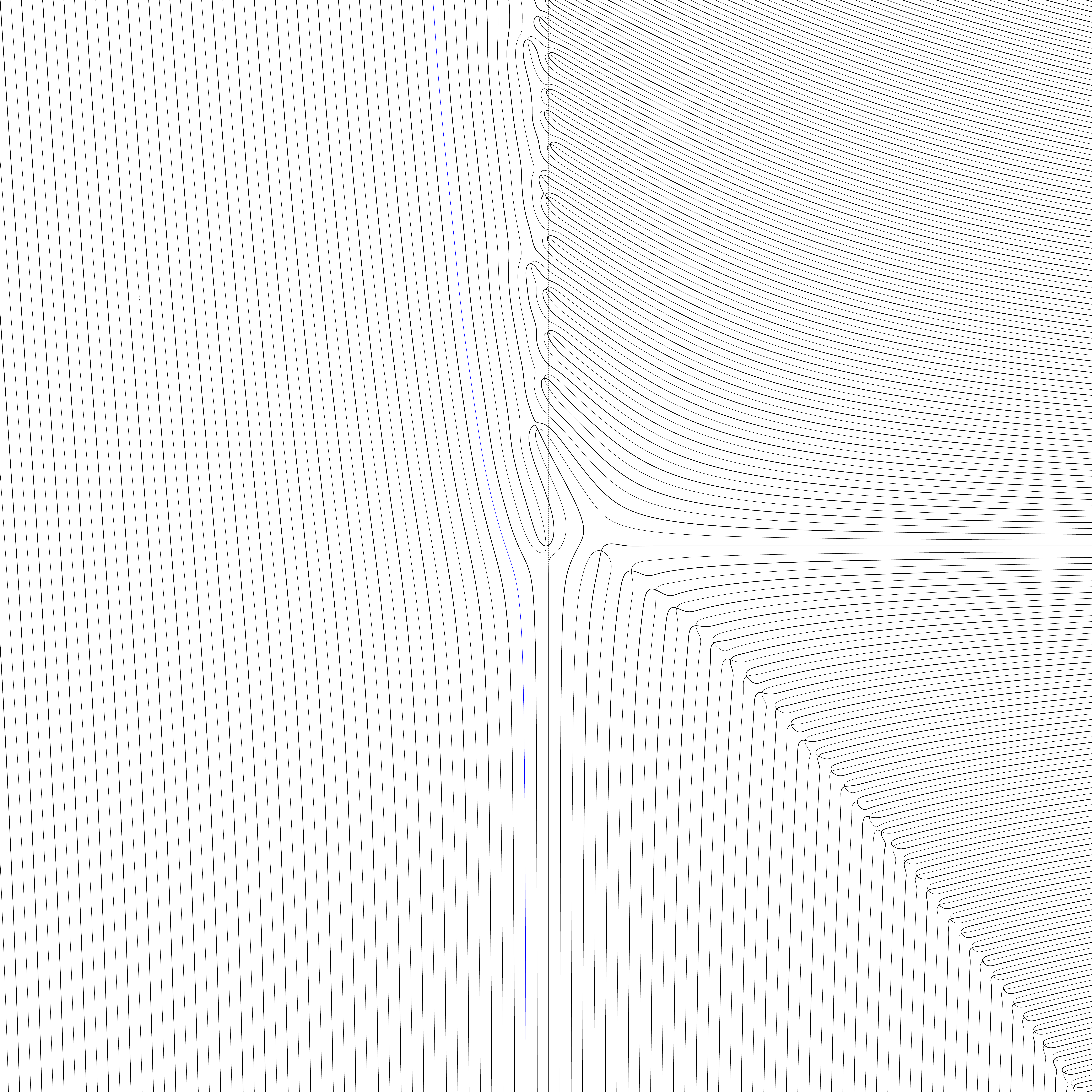}
\end{center}
\end{figure}

\begin{figure}[H]
\begin{center}
\includegraphics[width=\hsize]{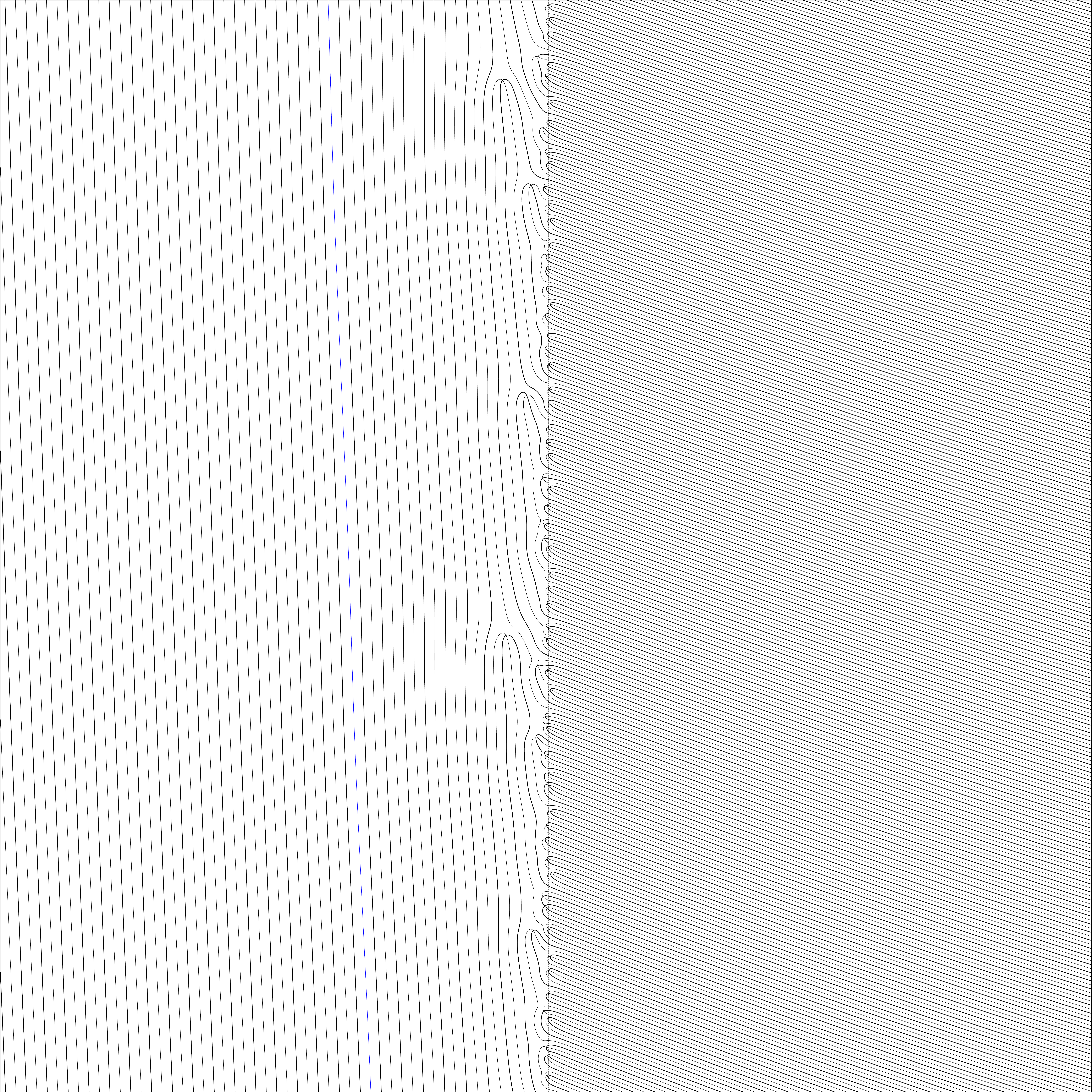}
\end{center}
\end{figure}

\end{document}